\documentclass[a4paper,11pt]{article}
\usepackage[nointlimits]{amsmath}
\usepackage{amssymb,amsmath,amsthm}
\usepackage{amsfonts}
\usepackage{latexsym}
\usepackage{graphicx}
\usepackage{color}
\usepackage{layout}
\usepackage[english]{babel}
\def\ds{\displaystyle}

\def\ZZ{\mathbb Z}

\def\bs{{\bf s}}

\newcommand{\leg}[2]{\left({#1\over #2}\right)}

\newtheorem{thm}{Theorem}[section]
\newtheorem{lem}[thm]{Lemma}

\newtheorem{conj}[thm]{Conjecture}

\begin{document}

\title{\bf Supercongruences\\ for a truncated hypergeometric series}
\author{{\sc Roberto Tauraso}\\
Dipartimento di Matematica\\
Universit\`a di Roma ``Tor Vergata'', Italy\\
{\tt tauraso@mat.uniroma2.it}\\
{\tt http://www.mat.uniroma2.it/$\sim$tauraso}
}
\date{}
\maketitle

\begin{abstract}
\noindent The purpose of this note is to obtain some congruences modulo
a power of a prime $p$ involving the truncated hypergeometric series 
$$\sum_{k=1}^{p-1} {(x)_k(1-x)_k\over (1)_k^2}\cdot{1\over k^a}$$
for $a=1$ and $a=2$. In the last section, the special case $x=1/2$ is considered.
\end{abstract}

\makeatletter{\renewcommand*{\@makefnmark}{}
\footnotetext{{2000 {\it Mathematics Subject Classification}: 33C20, 11B65, (Primary) 05A10, 05A19 (Secondary)}}\makeatother}

\section{Introduction}

In \cite{Mo:03}, E. Mortenson developed a general framework for studying
congruences modulo $p^2$ for the truncated hypergeometric series
$$F_{\!\!\!\!\!\!\!\!2\ \ 1}(x,1-x;1;1)_{\mbox{\small tr}(p)}
=\sum_{k=0}^{p-1} {(x)_k(1-x)_k\over (1)_k^2}.$$
Here we would like to present our investigations about a similar class of finite sums.
In order to state our main result we need to introduce the definition of what we call
{\sl Pochhammer quotient}
$$Q_p(x)={1\over p}\left(1-{(x)_p(1-x)_p\over (1)^2_p}\cdot {m^2\over [r/p]_m[-r/p]_m}\right)$$
where $p$ is an odd prime, $x=r/m$, $0<r<m$ are integers with $m$ prime to $p$,
$(x)_n=x(x+1)\cdots (x+n-1)$ denotes the Pochhammer symbol, and $[a]_m$ is the unique representative of $a$ modulo $m$ in $\{0,1,\dots,m-1\}$.
Note that the Pochhammer quotient is really a $p$-integral because by the partial-fraction decomposition
$${(1)_p\over (x)_p}=
\sum_{k=0}^{p-1} {\displaystyle\lim_{x\to -k}\left({(x+k)(1)_p\over (x)_p}\right)\over x+k}=
p\sum_{k=0}^{p-1} {p-1\choose k}{(-1)^k\over x+k}\equiv
p\sum_{k=0}^{p-1} {1\over x+k}\equiv {m\over [r/p]_m} \pmod{p}.$$
It is easy to verify that for some particular values of $x$
the quotient $Q_p(x)$ is connected with the usual {\sl Fermat quotient} $q_p(x)=(x^{p-1}-1)/p\,$:
\begin{align*}
&Q_p(1/2)={1\over p}\left(1-{1\over 16^p}{2p\choose p}^2\cdot
{4\over [1/p]_2[-1/p]_2}\right)\equiv -q_p(1/16)&\pmod{p^2},\\
&Q_p(1/3)={1\over p}\left(1-{1\over 27^p}{3p\choose p}{2p\choose p}\cdot
{9\over [1/p]_3 [-1/p]_3}\right)\equiv -q_p(1/27)&\pmod{p^2},\\
&Q_p(1/4)={1\over p}\left(1-{1\over 64^p}{4p\choose 2p}{2p\choose p}\cdot
{16\over [1/p]_4 [-1/p]_4}\right)\equiv -q_p(1/64)&\pmod{p^2},\\
&Q_p(1/6)={1\over p}\left(1-{1\over (16\cdot 27)^p}{6p\choose 2p}{4p\choose 2p}\cdot
{36\over [1/p]_6[-1/p]_6}\right)\equiv -q_p(1/(16\cdot 27))&\pmod{p^2}.
\end{align*}
where we used one of the equivalent statement of Wolstenholme's theorem,
namely, 
${ap\choose bp}\equiv {a\choose b}$ modulo $p^3$ for any prime $p>3$.

Our main goal in this paper is to show the following:
if $p>3$ is a prime then
\begin{align*}
&\sum_{k=1}^{p-1} {(x)_k(1-x)_k\over (1)_k^2}\cdot{1\over k}\equiv
Q_p(x) +{1\over 2}p\,Q_p(x)^2&\pmod{p^2},\\
&\sum_{k=1}^{p-1} {(x)_k(1-x)_k\over (1)_k^2}\cdot{1\over k^2}\equiv
-{1\over 2}Q_p(x)^2&\pmod{p}.
\end{align*}
The next section contains some preliminary results about certain hypergeometric identities and congruences.
We state and prove our main result in the third section. We conclude the paper by considering the
special case $x=1/2$ and by proving that
$$\sum_{k=1}^{p-1} {(1/2)_k^2\over (1)_k^2}\cdot{1\over k}\equiv-2H_{(p-1)/2}(1)\pmod{p^3}$$
where $H_n(r)=\sum_{k=1}^n {1\over k^r}$ is the {\sl finite harmonic sum} of order $n$ and weight $r$.

\section{Preliminaries}
The hypergeometric identity presented in the next theorem is known
(see, for example, equation~(21), Ch.~5.2 in~\cite{Luke:75}). 
Here we give a simple proof by using the Wilf-Zeilberger (WZ) method.

\begin{thm}\label{T1} For $0<x<1$ and for any positive integer $n$,
\begin{equation}\label{due}
\sum_{k=0}^{n-1} {(x)_k(1-x)_k\over (1)_k^2}\cdot{1\over n-k}={(x)_n(1-x)_n\over (1)_n^2}
 \left(\sum_{k=0}^{n-1}{1\over x+k}+\sum_{k=0}^{n-1}{1\over 1-x+k}\right).
\end{equation}
\end{thm}
\begin{proof} 
For $k=0\,\dots,n-1$, let
$$F(n,k)={(1)_n^2\over (x)_n(1-x)_n}\cdot{(x)_k(1-x)_k\over (1)_k^2}\cdot{1\over n-k}.$$
Then WZ method yields the certificate
$$R(n,k)=-{k^2(n-k)\over (n+1-k)(x+n)(1-x+n)}$$
with $G(n,k)=R(n,k)F(n,k)$, such that
$$F(n+1,k)-F(n,k)=G(n,k+1)-G(n,k).$$
It is easy to verify by induction that
$$\sum_{k=0}^{n-1}F(n,k)+\sum_{j=0}^{n-1}G(j,0)
=\sum_{k=0}^{n-1}\left(F(k+1,k)+G(k,k)\right).$$
Thus the desired identity follows by noting that $G(j,0)=0$ and
$$F(k+1,k)+G(k,k)={(1+k)^2\over (x+k)(1-x+k)}-{k^2\over (x+k)(1-x+k)}
={1\over x+k}+{1\over 1-x+k}.$$
\end{proof}

Let $B_n(x)$ be the {\sl Bernoulli polynomial} in $x$ of degree $n\geq 0$, given by
$$B_n(x)=\sum_{k=0}^n{n\choose k}B_k x^{n-k}$$
where $B_k$ are rational numbers called {\sl Bernoulli numbers} which are defined
recursively as follows:
$$B_0=1,\quad \mbox{and}\quad \sum_{k=0}^{n-1}{n\choose k}B_k=0\quad \mbox{for $n\geq 2$.}$$
For the main properties of $B_n(x)$ we will refer to Chapter 15 in \cite{IrRo:90}
and to the nice introductory article \cite{Ap:08}.

\begin{lem}\label{L1}
If $p>3$ is a prime and $x=r/m$, where $0<r<m$ are integers with $m$ prime to $p$, then for any positive integer $a$
\begin{align}\label{tre}
&\sum_{k=0}^{ap-1}\left(
{1\over x+k}+{1\over 1-x+k}\right)
-{m\over p}\sum_{j=0}^{a-1}\left(
{1\over [r/p]_m+jm}+{1\over [-r/p]_m+jm}\right)\nonumber\\
&\qquad\qquad\qquad\qquad\qquad\qquad
\equiv - {2\over 3}\,(ap)^2\,B_{p-3}(x)\pmod{p^3},\\
\label{quattro}
&{(x)_{ap}(1-x)_{ap}\over (1)^2_{ap}} 
\equiv {(x)_{(a-1)p}(1-x)_{(a-1)p}\over (1)^2_{(a-1)p}} \cdot {(x)_{p}(1-x)_{p}\over (1)^2_{p}} 
\cdot {1\over a^2}\left(1+{a(a-1)m^2\over [r/p]_m[-r/p]_m}\right)
\pmod{p^3}.
\end{align}
\end{lem}
\begin{proof}
We first show (\ref{tre}).
By well-known properties of the Bernoulli polynomials 
$$\sum_{k=0}^{ap-1}(x+k)^{\varphi(p^3)-1}
= {B_{\varphi(p^3)}(ap+x)-B_{\varphi(p^3)}(x)\over\varphi(p^3)}
={1\over\varphi(p^3)}\sum_{k=1}^{\varphi(p^3)}{\varphi(p^3)\choose k}B_{\varphi(p^3)-k}(x)(ap)^k.
$$
Since $m^nB_n(x)-B_n\in\ZZ$ (see \cite{Su:93} for a short proof),
it follows from the Clausen-von Staudt congruence (p. 233 in \cite{IrRo:90}) that
$$pB_n(x)\equiv {pB_n\over m^n}\equiv\left\{
\begin{array}{ll}
-1 &\mbox{if $(p-1)$ divides $n$ and $n>1$}\\
0 &\mbox{otherwise}
\end{array}
\right.\pmod{p}.$$
Moreover, by Kummer's congruences (p. 239 in \cite{IrRo:90}), if $q\geq 0$ and $2\leq r \leq p-1$ then
$${B_{q(p-1)+r}(x)\over q(p-1)+r}\equiv {B_r(x)\over r}\pmod{p}.$$
Hence, because $p$ divides the numerator of the fraction $x+k$ for $k=0,\dots,p-1$
if and only if $k=k_{0}$, we have that
\begin{align*}
\sum_{k=0}^{ap-1}{1\over x+k}-{m\over p}\sum_{j=0}^{a-1}{1\over [r/p]_m+jm}
&\equiv \sum_{k=0}^{ap-1}(x+k)^{\varphi(p^3)-1}\\
&\equiv ap\,B_{\varphi(p^3)-1}(x)- {1\over 2}\,(ap)^2\,B_{\varphi(p^3)-2}(x)\\
&\equiv ap\,B_{p^2(p-1)-1}(x)- {1\over 3}\,(ap)^2\,B_{p-3}(x)\pmod{p^3}.
\end{align*}
Finally, since by the symmetry relation $B_n(x)=(-1)^nB_n(1-x)$, we have that 
(\ref{tre}) holds.

As regards (\ref{quattro}), let $k_{0}:=p[r/p]_m/m-x$, then
\begin{align*}
{(x)_{ap}(1-x)_{ap}\over (x)_{(a-1)p}(1-x)_{(a-1)p}(x)_p(1-x)_p}
&=\prod_{k=0}^{p-1}\left(1+{(a-1)p\over x+k}\right)\left(1+{(a-1)p\over p-(x+k)}\right)\\
&=\prod_{k=0}^{p-1}\left(1+{a(a-1)p^2\over (x+k)(p-(x+k))}\right)\\
&\equiv 1+{a(a-1)p^2\over (x+k_0)(p-(x+k_0))}=1+{a(a-1)m^2\over [r/p]_m[-r/p]_m}\pmod{p^3}.
\end{align*}
Moreover, by an equivalent statement of Wolstenholme's theorem (see for example \cite{Zh:07}),
$${(1)_{ap}\over (1)_{(a-1)p}(1)_p}={ap\choose p}\equiv a \pmod{p^3}$$
and the proof of (\ref{quattro}) is complete.
\end{proof}

The next lemma provides a powerful tool which will be used in the next section
in the proof of the main theorem.

\begin{lem}\label{T2}
If $p>3$ is a prime and $x=r/m$, where $0<r<m$ are integers with $m$ prime to $p$, then
\begin{equation}\label{sei}
\sum_{k=ap+1}^{ap+p-1} {(x)_{k}(1-x)_{k}\over (1)_k^2}\cdot{1\over k}
\equiv {(x)_{ap}(1-x)_{ap}\over (1)^2_{ap}}
\sum_{k=1}^{p-1}{(x)_{k}(1-x)_{k}\over (1)_k^2}\cdot{1\over ap+k}
\pmod{p^2}.
\end{equation}
\end{lem}
\begin{proof}
Since
$${(x)_{ap+k}\over (x)_{ap}}=\prod_{j=0}^{k-1}(x+j+ap)\equiv (x)_k\biggl(1+ap\sum_{j=0}^{k-1}{1\over x+j}
\biggr)\pmod{p^2}$$
we have
$${(x)_{ap+k}(1-x)_{ap+k}(1)_{ap}^2
\over (x)_{ap}(1-x)_{ap}(1)_{ap+k}^2}\equiv
{(x)_{k}(1-x)_{k}\over (1)_{k}^2}
\left(1+ap\sum_{j=0}^{k-1}\left({1\over x+j}+{1\over 1-x+j}-{2\over 1+j}\right)\right)\pmod{p^2}.$$
Hence it suffices to prove that
$$ap\sum_{k=1}^{p-1}{(x)_{k}(1-x)_{k}\over (1)_{k}^2}\cdot{1\over ap+k}
\sum_{j=0}^{k-1}\left({1\over x+j}+{1\over 1-x+j}-{2\over 1+j}\right)\equiv 0\pmod{p^2},$$
that is, 
$$\sum_{k=1}^{p-1}{(x)_{k}(1-x)_{k}\over (1)_{k}^2}\cdot{1\over k}
\sum_{j=0}^{k-1}\left({1\over x+j}+{1\over 1-x+j}\right)\equiv
\sum_{k=1}^{p-1}{(x)_{k}(1-x)_{k}\over (1)_{k}^2}\cdot{2H_k(1)\over k}
\pmod{p}.$$
By (\ref{due}), the left-hand side becomes
\begin{eqnarray*}
\sum_{k=1}^{p-1}{1\over k}
\sum_{j=0}^{k-1}{(x)_{j}(1-x)_{j}\over (1)_{j}^2}\cdot{1\over k-j}
&=&
\sum_{j=0}^{p-2}{(x)_{j}(1-x)_{j}\over (1)_{j}^2}
\sum_{k=j+1}^{p-1}{1\over k(k-j)}\\
&=&
H_{p-1}(2)+\sum_{j=1}^{p-2}{(x)_{j}(1-x)_{j}\over (1)_{j}^2}
\left({1\over j}\sum_{k=j+1}^{p-1}\left({1\over k-j}-{1\over k}\right)\right)\\
&=&
H_{p-1}(2)+\sum_{j=1}^{p-2}{(x)_{j}(1-x)_{j}\over (1)_{j}^2}
\left({1\over j}\left(H_{p-1-j}(1)-H_{p-1}(1)+H_{j}(1)\right)\right)\\
&\equiv &
\sum_{j=1}^{p-2}{(x)_{j}(1-x)_{j}\over (1)_{j}^2}\cdot{2H_{j}(1)\over j}\equiv
\sum_{j=1}^{p-1}{(x)_{j}(1-x)_{j}\over (1)_{j}^2}\cdot{2H_{j}(1)\over j} \pmod{p},
\end{eqnarray*}
because $H_{p-1-j}(1)\equiv H_j(1)$ (mod $p$) and $H_{p-1}(1)\equiv H_{p-1}(2)\equiv 0$ (mod $p$).
The proof is now complete. 
\end{proof}

\section{The main theorem}\label{sec:main}

\begin{thm}\label{T3}
If $p>3$ is a prime and $x=r/m$, where $0<r<m$ are integers with $m$ prime to $p$, then
\begin{eqnarray}\label{A}
&&\sum_{k=1}^{p-1} {(x)_k(1-x)_k\over (1)_k^2}\cdot{1\over k}\equiv
Q_p(x) +{1\over 2}p\,Q_p(x)^2\pmod{p^2},\\\label{B}
&&\sum_{k=1}^{p-1} {(x)_k(1-x)_k\over (1)_k^2}\cdot{1\over k^2}\equiv
-{1\over 2}Q_p(x)^2\pmod{p}.
\end{eqnarray}
\end{thm}
\begin{proof} Let
$$S_a(b)=\sum_{k=ap+1}^{ap+p-1} {(x)_k(1-x)_k\over (1)_k^2}\cdot{1\over k^b}.$$
By (\ref{due}) with $n=p$ and by (\ref{tre}), we obtain
\begin{eqnarray}\label{sette}
-S_0(1)-pS_0(2)\equiv -{1\over p}+{(x)_{p}(1-x)_{p}\over (1)_{p}^2}\left(\sum_{k=0}^{p-1}{1\over x+k}+\sum_{k=0}^{p-1}{1\over 1-x+k}\right)
\equiv {\beta tm-1\over p}\pmod{p^2}.
\end{eqnarray}
where
$$\beta={(x)_{p}(1-x)_{p}\over (1)_{p}^2},\quad \mbox{and}\quad
t={1\over [r/p]_m}+{1\over [-r/p]_m}={m\over [r/p]_m [-r/p]_m}.$$
Moreover by (\ref{quattro})
$${(x)_{2p}(1-x)_{2p}\over (1)_{2p}^2}\equiv {\beta^2\over 4}\left(1+2mt\right)\pmod{p^3}.$$
and
$${1\over [r/p]_m+m}+{1\over [-r/p]_m+m}={3t\over 1+2mt}.$$
Hence, by (\ref{due}) with $n=2p$ and by (\ref{tre}), we get
\begin{eqnarray}
\lefteqn{-S_0(1)-2pS_0(2)-S_1(1)-2pS_1(2)}\nonumber\\
&&\ \ \ \ \
\equiv -{\beta\over p}-{1\over 2p}+{(x)_{2p}(1-x)_{2p}\over (1)_{2p}^2}
\left(\sum_{k=0}^{2p-1}{1\over x+k}+\sum_{k=0}^{2p-1}{1\over 1-x+k}\right)\nonumber\\
&&\ \ \ \ \ \label{otto}
\equiv -{\beta\over p}-{1\over 2p}+{\beta^2\over 4}\left(1+2mt\right)
{m\over p}\left(t+{3t\over 1+2mt}\right)\pmod{p^2}
\end{eqnarray}
According to Lemma~\ref{T2},
$$S_1(1)\equiv \beta(S_0(1)-pS_0(2)) \pmod{p^2},\quad S_1(2)\equiv \beta\,S_0(2)\pmod{p},$$
and (\ref{sette}) and (\ref{otto}) yield the
linear system
$$\left\{
\begin{array}{l}
\ds-S_0(1)-pS_0(2)\equiv {\beta tm-1\over p}\pmod{p^2},\\\\
\ds-S_0(1)-2pS_0(2)-\beta S_0(1)-\beta\,pS_0(2)
\equiv-{\beta\over p}-{1\over 2p}+{\beta^2\over 4}(1+2mt){m\over p}
\left(t+{3t\over 1+2mt}\right)\pmod{p^2}.
\end{array}
\right.$$
By substituting the first congruence in the second one, we obtain
$$
-pS_0(2)+{\beta tm-1\over p}+
\beta{\beta tm-1\over p}
\equiv -{\beta \over p}-{1\over 2p}+{\beta^2\over 4}(1+2mt){m\over p}
\left(t+{3t\over 1+2mt}\right)
\pmod{p^2},$$
which yields
$$S_0(2)\equiv -{1\over 2}\left({1-\beta tm\over p}\right)^2 \equiv -{1\over 2}Q_p(x)^2\pmod{p}.$$
Finally,
$$S_0(1)\equiv Q_p(x) +{1\over 2}p\,Q_p(x)^2\pmod{p^2}.$$
\end{proof}

We conclude this section by posing a conjecture which extends (\ref{B}).

\begin{conj}
If $p>3$ is a prime and $x=r/m$, where $0<r<m$ are integers with $m$ prime to $p$, then
$$\sum_{k=1}^{p-1} {(x)_k(1-x)_k\over (1)_k^2}\cdot{1\over k^2}\equiv
-{1\over 2}Q_p(x)^2-{1\over 2}p\,Q_p(x)^3\pmod{p^2}.$$
\end{conj}

More conjectures of the same flavour can be found in Section~A30 of~\cite{Sunzw:10}.

\section{The special case $x={1\over 2}$}\label{sec:special}

As we already noted, if $x={1\over 2}$ then 
$${(x)_k(1-x)_k\over (1)_k^2}={(1/2)_k^2\over (1)_k^2}={2k\choose k}^2{1\over 16^k}.$$
Let $p$ and odd prime then for $0\leq k \leq n=(p-1)/2$ we have that
$${n+k\choose 2k}={\prod_{j=1}^k(p^2-(2j-1)^2)\over 4^k(2k)!}\equiv
{\prod_{j=1}^k(2j-1)^2\over (-4)^k(2k)!}={2k\choose k}{(-1)^k\over 16^k}\pmod{p^2},$$
which means that
$$(-1)^k{n\choose k}{n+k\choose k}=(-1)^k{2k\choose k}{n+k\choose 2k}
\equiv{2k\choose k}^2{1\over 16^k}\pmod{p^2}.$$
Since $p$ divides ${2k\choose k}$ for $n<k<p$, it follows that for any $p$-adic integers $a_0,a_1,\dots a_{p-1}$
$$\sum_{k=0}^{p-1}{2k\choose k}^2{a_k\over 16^k}\equiv \sum_{k=0}^{n}(-1)^ka_k{n\choose k}{n+k\choose k} \pmod{p^2}.$$
This remark is interesting because the sum on the right-hand side could be easier to study modulo $p^2$.
With this purpose in mind, we consider the identity 2.1 in \cite{Pr:07}
\begin{equation}\label{auno}
\sum_{k=1}^n {(-1)^k\over z+k}{n\choose k}{n+k\choose k}=
{1\over z}\left({(1-z)_n\over (1+z)_n}-1\right).
\end{equation}
As a first example, we can give a short proof of (1.1) in \cite{Mo:03}.
By multiplying (\ref{auno}) by $z$, and by letting $z\to\infty$, we obtain that  
$$\sum_{k=0}^n {(-1)^k}{n\choose k}{n+k\choose k}=(-1)^n$$
which implies that for any prime $p>3$
$$\sum_{k=0}^{p-1}{2k\choose k}^2{1\over 16^k}\equiv (-1)^{p-1\over 2}=\leg{-1}{p} \pmod{p^2}.$$
For more applications of the above remark see \cite{Sunzh:10} and \cite{Lo:10}.

Let $\bs=(s_1,s_2,\dots,s_d)$ be a vector whose entries are positive integers then we
define the {\sl multiple harmonic sum} for $n\geq 0$ as
$$H_n(s_1,s_2,\dots,s_d)=\sum_{1\leq k_1<k_2<\cdots<k_d\leq n} \!{1\over k_1^{s_1}k_2^{s_2}\cdots k_d^{s_d}}.$$
We call $d$ and $|\bs|=\sum_{i=1}^d s_i$ its depth and its weight respectively.

\begin{thm}\label{T4} For $n,r\geq 1$
\begin{equation}\label{adue} 
\sum_{k=1}^n {(-1)^k\over k^r}{n\choose k}{n+k\choose k}=
-\sum_{d=1}^r 2^d\sum_{|\bs|=r} H_n(s_1,s_2,\dots,s_d).
\end{equation}
\end{thm}
\begin{proof} 
Note that
$$(1+z)_n=n!\,\left(1+\sum_{d=1}^n H_n(\overbrace{1,1,\dots,1}^{\mbox{$d$ times}})z^d\right),\quad\mbox{and}\quad
{d\over dz}(1+z)_n=(1+z)_n\sum_{k=0}^{n-1}{1\over z+k}.$$
Then the left-hand side of (\ref{adue}) can be obtained by differentiating $(r-1)$ times 
the identity (\ref{auno}) with respect to $z$ and then by letting $z=0$:
\begin{align*}
\sum_{k=1}^n {(-1)^k\over k^r}{n\choose k}{n+k\choose k}&=
{(-1)^{r-1}\over (r-1)!}\left({d^{r-1}\over dz^{r-1}}\left({1\over z}\left({(1-z)_n\over (1+z)_n}-1\right)\right)\right)_{z=0}\\
&={(-1)^{r-1}\over r!}\left({d^{r}\over dz^{r}}\left({(1-z)_n\over (1+z)_n}\right)\right)_{z=0}.
\end{align*}
Taking the derivatives we get a formula which involves products of multiple harmonic sums. 
This formula can be simplified to the right-hand side of (\ref{adue}) by using the so-called {\sl stuffle} product (see for example \cite{Ta:10b}).
\end{proof}

In the next theorem, we prove two generalizations of (\ref{A}) and (\ref{B}) for $x=1/2$. 
Note that the first of these congruences can be considered as a variation of another congruence proved by the author
in \cite{Ta:10a}: for any prime $p>3$
$$\sum_{k=1}^{p-1} {(1/2)_k\over (1)_k}\cdot{1\over k}\equiv
-H_{(p-1)/2}(1)\pmod{p^3}.$$

\begin{thm}\label{T5} For any prime $p>3$
\begin{align*}
&\sum_{k=1}^{p-1} {(1/2)_k^2\over (1)_k^2}\cdot{1\over k}\equiv
-2H_{(p-1)/2}(1)&\pmod{p^3},\\
&\sum_{k=1}^{p-1} {(1/2)_k^2\over (1)_k^2}\cdot{1\over k^2}\equiv
-2H_{(p-1)/2}(1)^2&\pmod{p^2}.
\end{align*}
\end{thm}

\begin{proof} We will use the same notations as in Theorem~\ref{T3}. For $x=1/2$ we have that $m=t=2$. 
By (\ref{due}) with $n=p$ and by (\ref{tre}), we obtain
\begin{align}-S_0(1)-pS_0(2)-p^2S_0(3)
&\equiv -{1\over p}+\beta\left(\sum_{k=0}^{p-1}{1\over x+k}+\sum_{k=0}^{p-1}{1\over 1-x+k}\right)\nonumber\\
&\equiv {\beta tm-1\over p}-{2\beta\over 3}\,p^2\,B_{p-3}(1/2)\pmod{p^3}.\label{atre}
\end{align}
On the other hand, by \cite{Zh:07}
\begin{align*}
\beta&\equiv{1\over 16^p}{2p\choose p}^2\equiv {4\over 16^p}\left(1-{2\over 3}\,p^3\, B_{p-3}\right)^2
\equiv {1\over 4(1+pq_p(2))^4}\left(1-{4\over 3}\,p^3\, B_{p-3}\right)
\pmod{p^4},
\end{align*}
and by Raabe's multiplication formula
$$B_{p-3}(1/2)=\left(1-{1\over 2^{p-4}}\right) B_{p-3}\equiv 7\,B_{p-3}\pmod{p}.$$
Moreover, by letting $n=(p-1)/2$ in (\ref{adue}) for $r=2,3$, we have that
\begin{align*}
S_0(2)&\equiv\sum_{k=1}^n {(-1)^k\over k^2}{n\choose k}{n+k\choose k}=
-2H_n(2)-4H_n(1,1)=-2\,H_n(1)^2\pmod{p^2},\\
S_0(3)&\equiv\sum_{k=1}^n {(-1)^k\over k^3}{n\choose k}{n+k\choose k}=
-2H_n(3)-4H_n(2,1)-4H_n(1,2)-8H_n(1,1,1)\\
&=-{4\over 3}\,H_n(1)^3-{2\over 3}\,H_n(3)\pmod{p^2}.
\end{align*}
Finally by \cite{Sunzh:00}
$$H_{n}(1)\equiv -2q_p(2)-pq_p(2)^2+{2\over 3}\,p^2q_p(2)^3+{7\over 12}\,p^2B_{p-3}\pmod{p^3}\,,
\quad H_{n}(3)\equiv 2B_{p-3} \pmod{p}.$$
By plugging all of these values in (\ref{atre}), after a little manipulation, we easily verify the desired congruence
for $S_0(1)$.
\end{proof}


\begin{thebibliography}{99}

\bibitem{Ap:08} T. M. Apostol,
{\it A primer on Bernoulli numbers and polynomials},
Math. Mag. {\bf 81} (2008), 178--190.

\bibitem{IrRo:90} K. Ireland, M. Rosen,
{\it A Classical Introduction to Modern Number Theory}, Springer, New York, 1990.

\bibitem{Lo:10} L. Long, 
{\it Hypergeometric evaluation identities and supercongruences},
preprint {\sf arXiv:math.NT/arXiv:0912.0197} (2010).

\bibitem{Luke:75} Y. L. Luke,
{\it Mathematical Functions and Their Approximations}, Academic Press,
New York, 1975.

\bibitem{Mo:03}  E. Mortenson,
{\it Supercongruences between truncated
$\ F_{\!\!\!\!\!\!\!\!2\ \ 1}$ hypergeometric functions and
their Gaussian analogs}, Trans. Amer. Math. Soc.  {\bf 355} (2003), 987--1007.

\bibitem{Pr:07}  H. Prodinger,
{\it Human proofs of identities by Osburn and Schneider}, 
preprint {\sf arXiv:math.CO/0710.0464} (2007).

\bibitem{Sunzh:00} Z. H. Sun, 
{\it Congruences concerning Bernoulli numbers and Bernoulli polynomials}, 
Discrete Appl. Math. {\bf 105} (2000), 193--223. 

\bibitem{Sunzh:10} Z. H. Sun, 
{\it Congruences concerning Legendre polynomials},
Proc. Amer. Math. Soc. {\bf 139} (2011) 1915--1929. 

\bibitem{Sunzw:10} Z. W. Sun, 
{\it Open conjectures on congruences},
preprint {\sf arXiv:math.NT/0911.5665} (2011).

\bibitem{Su:93} B. Sury, 
{\it The value of Bernoulli polynomials at rational numbers},
Bull. London Math. Soc.  {\bf 25} (1993), 327--329.

\bibitem{Ta:10a} R. Tauraso,
{\it Congruences involving alternating multiple harmonic sum},
Electron. J. Combin., R16 (2010).

\bibitem{Ta:10b} R. Tauraso,
{\it New harmonic number identities with applications},
S\'em. Lothar. Combin., B63g (2010).

\bibitem{Zh:07} J. Zhao,
{\it Bernoulli numbers, Wolstenholme's theorem, and $p^5$ variations of Lucas' theorem},
J. Number Theory {\bf 123} (2007), 18--26.

\end{thebibliography}
\end{document}